\documentstyle[11pt]{article}

\input amssym.def
\input amssym

\begin{document}

\def \1{{\bf 1}}
\def \a{{{\frak a}}}
\def \al{\alpha}
\def \ar{{\alpha_r}}
\def \A{{\Bbb A}}
\def \Ad{{\rm Ad}}
\def \b{{{\frak b}}}
\def \bs{\backslash}
\def \B{{\cal B}}
\def \cent{{\rm cent}}
\def \C{{\Bbb C}}
\def \CA{{\cal A}}
\def \CB{{\cal B}}
\def \CE{{\cal E}}
\def \CF{{\cal F}}
\def \CG{{\cal G}}
\def \CH{{\cal H}}
\def \CL{{\cal L}}
\def \CM{{\cal M}}
\def \CN{{\cal N}}
\def \CP{{\cal P}}
\def \CQ{{\cal Q}}
\def \CO{{\cal O}}
\def \CS{{\cal S}}
\def \det{{\rm det}}
\def \diag{{\rm diag}}
\def \dom{{\rm dom}\hspace{2pt}}
\def \e{\epsilon}
\def \End{{\rm End}}
\def \Fx{{\frak x}}
\def \FX{{\frak X}}
\def \g{{{\frak g}}}
\def \ga{\gamma}
\def \Ga{\Gamma}
\def \GL{{\rm GL}}
\def \h{{{\frak h}}}
\def \Hom{{\rm Hom}}
\def \im{{\rm im}\hspace{2pt}}
\def \Ind{{\rm Ind}}
\def \k{{{\frak k}}}
\def \K{{\cal K}}
\def \la{\lambda}
\def \lap{\triangle}
\def \Lie{{\rm Lie}\hspace{2pt}}
\def \m{{{\frak m}}}
\def \mathqed{\hfill \Box} 
\def \mod{{\rm mod}}
\def \Mat{{\rm Mat}}
\def \n{{{\frak n}}}
\def \name{\bf}
\def \N{\Bbb N}
\def \ord{{\rm ord}}
\def \O{{\cal O}}
\def \p{{{\frak p}}}
\def \ph{\varphi}
\def \prf{{\bf Proof: }}
\def \qed{\hfill {$\Box$} 

$ $

}
\def \Q{\Bbb Q}
\def \res{{\rm res}}
\def \R{{\Bbb R}}
\def \Re{{\rm Re \hspace{1pt}}}
\def \ra{\rightarrow}
\def \rank{{\rm rank}}
\def \Rep{{\rm Rep}}
\def \sign{{\rm sign}\hspace{2pt}}
\def \supp{{\rm supp}}
\def \t{{{\frak t}}}
\def \T{{\Bbb T}}
\def \tr{{\hspace{1pt}\rm tr\hspace{1pt}}}
\def \vol{{\rm vol}\hspace{1pt}}
\def \V{{\cal V}}
\def \z{\frak z}
\def \Z{\Bbb Z}
\def \={\ =\ }

\newcommand{\rez}[1]{\frac{1}{#1}}
\newcommand{\der}[1]{\frac{\partial}{\partial #1}}
\newcommand{\binom}[2]{\left( \begin{array}{c}#1\\#2\end{array}\right)}
\newcommand{\norm}[1]{\hspace{2pt}\parallel #1 \parallel\hspace{2pt}}

\newcounter{lemma}
\newcounter{corollary}
\newcounter{proposition}
\newcounter{theorem}
\newcounter{zwisch}

\renewcommand{\subsection}{\refstepcounter{subsection}\stepcounter{lemma} 
	\stepcounter{corollary} \stepcounter{proposition}
	\stepcounter{conjecture}\stepcounter{theorem}
	$ $ \newline
	{\bf \arabic{section}.\arabic{subsection}\hspace{8pt}}}

\newtheorem{conjecture}{\stepcounter{lemma} \stepcounter{corollary} 	
	\stepcounter{proposition}\stepcounter{theorem}
	\stepcounter{subsection}\hskip-12pt Conjecture}[section]
\newtheorem{lemma}{\stepcounter{conjecture}\stepcounter{corollary}	
	\stepcounter{proposition}\stepcounter{theorem}
	\stepcounter{subsection}\hskip-7pt Lemma}[section]
\newtheorem{corollary}{\stepcounter{conjecture}\stepcounter{lemma}
	\stepcounter{proposition}\stepcounter{theorem}
	\stepcounter{subsection}\hskip-7pt Corollary}[section]
\newtheorem{proposition}{\stepcounter{conjecture}\stepcounter{lemma}
	\stepcounter{corollary}\stepcounter{theorem}
	\stepcounter{subsection}\hskip-7pt Proposition}[section]
\newtheorem{theorem}{\stepcounter{conjecture} \stepcounter{lemma}
	\stepcounter{corollary}\stepcounter{proposition}		
	\stepcounter{subsection}\hskip-11pt Theorem}[section]

\title{A Polya operator for automorphic L-functions}
\author{{\small by}\\ {} \\ Anton Deitmar}
\date{}
\maketitle

\pagestyle{myheadings}
\markright{POLYA OPERATOR FOR AUTOMORPHIC L-FUNCTIONS}

\tableofcontents

\begin{center}
{\bf Introduction}
\end{center}

$$ $$

\section{Schwartz-Bruhat functions and Haar measures}
Let $k$ be a number field and let $M$ be a simple $k$-algebra with center $k$ of rank $n$, i.e. $\dim_k(M)=n^2$.
Fix a $k$-involution $x\mapsto x^\tau$ on $M$, i.e. $(xy)^\tau =y^\tau x^\tau$ as well as $(x^\tau)^\tau =x$ and $x^\tau =x$ when $x\in k$.
We will assume that $M$ and $\tau$ are defined over the ring of integers $\CO$ of $k$.
The standard example will be the algebra of $n\times n$ matrices $\Mat_n(k)$ and $x^\tau =x^t$ the transposed matrix.

For any ring $R$ over $k$ let $M(R)=M\otimes_kR$ and let $G$ be the $k$-groupscheme of invertible elements, i.e. $G(R) =M(R)^\times$.
For almost all places $v$ of $k$ we have $M(k_v)\cong \Mat_n(k_v)$.
In this case we say that $M$ {\bf splits at} $v$.
At any splitting place $v$ we assume a $\CO_v$-isomorphism $M(k_v)\ra \Mat_n(k_v)$ fixed, wher $\CO_v$ is the local ring of integers.

Let $\A =\A_f\times \A_\infty$ be the ring of adeles of $k$ where $\A_f$ is the finite and $\A_\infty$ the infinite adeles.
Let $dx$ denote the additive Haar-measure on $\A$ given by $dx=\otimes_{v}dx_v$, the product being extended over all places of $k$, where at a finite place $v$ we have for the valuation ring $\CO_v\subset k_v$ that $dx_v(\CO_v)=1$.
At the infinite places we will normalize the measures in a way such that the lattice $k\subset \A$ has covolume $1$, i.e. $\vol(\A /k)=1$.
We also write $dx=\otimes_vdx_v$ for an additive Haar-measure on $M(\A)$, where $dx_v$ is given componentwise by the above at each splitting place $v$.

For $a\in\A$ let $|a|$ be its modulus, so $dx(aA)=|a|dx(A)$.
For $x\in M(\A)$ let $|x|=|\det(x)|$, where $\det : M\ra k$ is the reduced norm. Note that $\det$ eqals the determinant over each field splitting $M$ and that $\det(x^\tau)=\det(x)$ for any $x$.
Let $dx^\times=\otimes_vdx^\times_v$ be the Haar measure on $G(\A)$ given by $dx_v^\times(G(\CO_v))=1$ for $v$ finite and $dx_v^\times=\frac{dx_v}{|x|_v^n}$ at the infinite places.
The measure
$$
|x|^ndx^\times \= \bigotimes_{v< \infty}|x|_v^n dx_v^\times \otimes \bigotimes_{v|\infty}dx_v
$$
is translation invariant on $M(\A)$ but it is not a Haar measure since it is infinity on any open set.

Let $N(\A)$ be the set of all $m\in M(\A)$ with $\det(m) =0$.
Then for $y\in M(k)$ and $m\notin N(\A)$ we have $ym\ne m$.
Let $\CS(M(\A))$ be the space of Schwartz-Bruhat functions on $M(\A)$, that is, any $f\in\CS(N(\A))$ is a finite sum of functions of the form $f=\otimes_vf_v$, where $f_v$ is the characteristic function of the set $\CO_v$ for almost all $v$ and $f_v\in\CS(M(k_v))$ at all places, where $\CS(M(k_v))$ is the usual Schwartz-Bruhat space if $v$ is infinite and is the space of locally constant functions of compact support if $v$ is finite.

To define Fourier transforms we will fix an additive character $\psi$ as follows.
At first assume $k=\Q$, then $\psi=(\prod_p\psi_p)\psi_\infty$ with
$\psi_p(\Z_p)=1$, $\psi_p(p^{-n})=e^{2\pi i/p^n},$
and $\psi_\infty(x)=e^{2\pi ix}.$
For general $k$ note that the trace map $Tr_{k/\Q}:k\ra \Q$ induces a trace $Tr: \A_k\ra \A_\Q$ and let $\psi_k :=\psi_\Q \circ Tr$.
The character $\psi$ identifies $\A$ with its dual via the pairing $\langle x,y\rangle =\psi(xy)$.
Note that $\psi$ is chosen in a way that the lattice $k\subset\A$ is its own dual, i.e. 
$$
\langle x,y\rangle\=1\ \ \forall y\in k\hspace{20pt} \Leftrightarrow\hspace{20pt} x\in k.
$$
For $f\in\CS(\A)$ its Fourier transform is defined by
$$
\hat{f}(x) \= \int_\A f(y)\psi(xy)dy.
$$
We lift these notions to $M(\A)$.
Let $\psi_M : M(\A)\ra\C$ be defined by $\psi_M(x) =\psi(\tr_{M/k}(x))$.
Then $\psi_M$ sets $M(\A)$ in self duality and $M(k)$ is a self dual lattice.
The Fourier transform for $f\in\CS(M(\A))$ is
$$
\hat{f} \= \int_{M(\A)}f(y)\psi_M(xy)dy.
$$
Let $\CS_0 =\CS(M(\A))_0$ be the space of all $f\in\CS(M(\A))$ such that $f$ and $\hat{f}$ send $N(\A)$ to zero.

\section{Local factors}\label{local}
Let $K$ be a nonarchimedean local field with ring of integers $\CO$.
We assume $K\supset k$ and $M$ splits over $K$.
So in this section read $M(K)=\Mat_n(K)$ and $G(K)=\GL_n(K)$.
Let $\varpi$ be a uniformizing element and let $|.| : K\ra \R_{\ge 0}$ be the absolute value normalized so that $|\varpi|=q^{-1}$, where $q$ is the number of elements of the residue class field $\CO /\varpi\CO$.

Let $\pi$ be an irreducible admissible Hilbert representation of $G(K)$, that is, $\pi$ is a continuous representation of $G(K)$ on a Hilbert space $V_\pi$ with $\dim V_\pi^U <\infty$ for any open subgroup $U$ of $G(K)$, where $V^U$ is the set of vectors fixed by any element of $U$.

Suppose $\pi$ is a {\bf class one representation}, i.e. $\pi$ has nonzero fixed vectors under the maximal compact subgroup $G(\CO)$ of $G(K)$.
The {\bf Hecke algebra} $\CH =\CH(G(K),G(\CO))$ is the convolution algebra of all compactly supported functions $f :G(K)\ra \C$ which are bi-invariant under $G(\CO)$, that is $f(xk)=f(kx)=f(x)$ for all $k\in G(\CO)$.

Let $A\subset G(K)$ be the subgroup of diagonal elements and $N\subset G(K)$ be the subgroup of upper triangular matrices with ones on the diagonal.
Then $B=AN$ is the Borel subgroup of upper triangular matrices.
For $a=\diag(a_1,\dots ,a_n)\in A$ let 
$$
\delta(a) := |\det (a|\Lie N)| = |a_1|^{n-1} |a_2|^{n-3} \dots |a_n|^{-(n-1)}
$$
the modulus of $B$.
Let $dn$ be the Haar measure on $N$ normalized so that $\vol(N\cap G(\CO))=1$.
The {\bf Satake transform} (see \cite{cart}, p. 146)
$$
Sf(a) \ :=\ \delta(a)^{\rez{2}}\int_Nf(an)dn
$$
gives an isomorphism $\CH \ra C[\lambda]^W$, where $\lambda = A/A\cap G(\CO)$ and $\C[\lambda]$ is the convolution algebra (group algebra) of all finitely supported functions on $\lambda$.
Finally $W$ is the Weyl group of all permutations of the diagonal entries of $A$.
The set of characters $\Hom (\lambda ,\C^*)$ forms  the complex points of a torus $T\subset \GL_n(\C)$.
The Hecke algebra $\CH$ acts on the one dimensional space $V_\pi^{G(\CO)}$ through a character $\chi_\pi$ which by the Satake isomorphism is given by an element of $T$.
The {\bf local L-factor} of $\pi$ is defined by
$$
L(\pi) \= \det(1-\chi_\pi)^{-1}.
$$
We will need to make this more explicit.
Let $\varpi_j = \diag(1,\dots ,\varpi,\dots ,1)$ (the $\varpi$ on the $j$-th position.
The Satake isomorphism gives a bijection
$$
\Hom_{alg}(\CH ,\C) \ra \Hom(\lambda ,\C^*)/W
$$
and using these terms the local factor is given by
\begin{eqnarray*}
L(\pi)^{-1} &=& \det(1-\chi_\pi)\\
	&=& \prod_j (1-\chi_\pi(\varpi_j)).
\end{eqnarray*}
Let $e\in\C[\lambda]^W$ be given by
$e=\prod_j (1-\varpi_j)$.
Then there exists a unique $f\in\CH$ with $Sf=e$ and for this $f$ we have
$$
\tr\pi(f) \= L(\pi)^{-1}.
$$

But also $L(\pi)$ itself rather then its inverse can sometimes be given as a trace.
Let $|\pi|:= \max_j(\chi_\pi(\varpi_j)).$
The characteristic function $\1_{M(\CO)}$ restricted to $G(K)$ does not lie in $\CH$ since it is not compactly supported but it can be written as an infinite sum of elements of $\CH$.
The same applies to the function $\1_{M(\CO)}|x|^{\frac{n-1}{2}}$.

\begin{proposition}
If $|\pi|<1$ then $\pi(\1_{M(\CO)}|x|^{\frac{n-1}{2}})$ exists and is of trace class with
$$
\tr\pi(\1_{M(\CO)}|x|^{\frac{n-1}{2}}) \= L(\pi).
$$
\end{proposition}

\prf
We compute the Satake transform of $f(x)=\1_{M(\CO)}(x)|x|^{\frac{n-1}{2}}$.
At first note that $f(an)\ne 0$ implies $a\in A\cap M(\CO)$, that is $|a_j|\le 1$ for $j=1,\dots ,n$.
For such an $a$ we substitute $y_k=a_jn_{j,k}$, where $n_{j,k}$ is the corresponding entry of $n$.
Since $\int_Nf(n)dn=1$ this gives
$$
\int_Nf(an)dn\= |a|^{\frac{n-1}{2}}|a_1|^{-(n-1)}|a_2|^{-(n-2)}\dots |a_{n-1}|^{-1}.
$$
And so
\begin{eqnarray*}
Sf(a) &=& \delta(a)^{\rez{2}}\int_Nf(an)dn\\
	&=& |a|^{\frac{n-1}{2}}|a_1 a_2 \dots a_n|^{-\frac{n-1}{2}}\\
	&=& 1.
\end{eqnarray*}
So that $Sf = \1_{A\cap M(\CO)}$.

The condition $|\pi|<1$ implies that
\begin{eqnarray*}
L(\pi) &=& \det(1-\chi_\pi)^{-1}\\
	&=& \prod_j(1-\chi_\pi(\varpi_j))^{-1}\\
	&=& \prod_j\sum_{k=0}^\infty \chi_\pi(\varpi_j)^k\\
	&=& \chi_\pi\left( \prod_j\sum_{k=0}^\infty \varpi_j^k\right)\\
	&=& \chi_\pi(\1_{A\cap M(\CO)}).
\end{eqnarray*}
\qed

The unramified character $g\mapsto |g|^s$ for some $s\in \C$ is a class one admissible representation of $G$ and so is $\pi_s=|.|^s\pi : g\mapsto |g|^s\pi(g)$.
We compute
$$
\chi_{\pi_s}\= q^{-s}\chi_\pi
$$
and
$$
\left|\pi_s\right| = q^{-\Re(s)} |\pi|,
$$
so that we have the

\begin{corollary}
For any irreducible admissible $\pi$ and $s\in \C$ with $|\pi|<q^{\Re(s)}$ it holds
\begin{eqnarray*}
L(\pi,s) &:=& L(\pi_s)\\
	 &=& \tr |.|^s\pi\left(\1_{M(\CO)}(x)|x|^{\frac{n-1}{2}}\right)\\
	 &=& \tr \pi\left(\1_{M(\CO)}(x)|x|^{s+\frac{n-1}{2}}\right).
\end{eqnarray*}
\end{corollary}
\qed

\section{Global $L$-functions}
For $f\in\CS(M(\A))_0$ let $f^\tau(x) := f(x^\tau)$.
Let then $E(f) :G(\A)\ra\C$ be defined by
\begin{eqnarray*}
E(f) &=& |x|^{\frac{n}{2}}\sum_{\ga\in G(k)}f(\ga x)\\
	&=& |x|^{\frac{n}{2}}\sum_{\ga\in M(k)}f(\ga x).
\end{eqnarray*}
Note that $f\in\CS(M(\A))_0$ implies that $f$ vanishes on the set $(M(k)-G(k))M(\A)$. 
Therefore it doesn't matter whether the sum is extended over $G(k)$ or $M(k)$.

\begin{lemma}
For any $f\in\CS(M(\A))_0$ the sum $E(f)(x)$
converges locally uniformly in $x$ with all derivatives. 
For $x\in G(\A)$ it holds
$$
E({f})(x)\= E(\hat{f}^{\tau})(x^{-\tau}),
$$
where $x^{-\tau}=(x^{-1})^\tau=(x^\tau)^{-1}$.
Finally for any $n\in\N$ there is a $C>0$ such that
$$
|E(f)(x)|\ \le\ C\min(|x|,\rez{|x|})^n.
$$
\end{lemma}

\prf
With $f$ any derivative of $f$ is again in $\CS(M(\A))$ so we only need to check convergence for $f$ itself.
Now $M(k)$ is a lattice in $M(\A)$and so is $M(k)x$ for any $x\in G(\A)$.
This lattice further depends continuously on $x$, which gives the claim.

For any lattice $\Ga\subset M(\A)$ and any $f\in\CS(M(\A))$ the Poisson summation formula says
$$
{\rm covol}(\Ga) \sum_{\ga\in\Ga} f(\ga) \= \sum_{\eta\in \Ga^\perp}\hat{f}(\eta),
$$
where $\Ga^\perp$ is the dual lattice.
Let $x\in G(\A)$, then $M(k)x$ still is a lattice in $M(\A)$ and its dual is $x^{-1}M(k)$.
Further, the covolume of $M(k)x$ is $|x|^n$ so that
\begin{eqnarray*}
|x|^n\sum_{\ga\in M(k)}f(\ga x) &=& \sum_{\Ga\in M(k)}\hat{f}(x^{-1}\ga)\\
	&=& \sum_{\Ga\in M(k)}\hat{f}^\tau(\ga x^{-\tau}).
\end{eqnarray*}

Now it remains to control the growth of $E(f)$ and by the last assertion it suffices to control the growth when $|x|$ is large.
We may assume that $f$ is of the form $f=f_{fin}\otimes f_\infty$ with $f_{fin}\in\CS(M(\A_{fin}))$ and $f_\infty\in\CS(M(\A_\infty))$.
The function $f_{fin}$ has compact support, say $C\subset M(\A_{fin})$.
Shifting factors we may assume that $|f_{fin}|\le \1_C$ and so we can estimate
\begin{eqnarray*}
|E(f)(x)|&\le& |x|^{\frac{n}{2}}\sum_{\ga\in G(k)} |f(\ga x)|\\
	&\le& |x|^{\frac{n}{2}}\sum_{\ga\in G(k)\cap Cx_{fin}^{-1}} |f_\infty(\ga x_\infty)|.
\end{eqnarray*}
There is $\tilde{x}_{fin}\in G(k)$ such that $Cx_{fin}^{-1}=c\tilde{x}_{fin}^{-1}$ and $|x_{fin}|_{fin}=|\tilde{x}_{fin}|_{fin}$ and so $|E(f)(x)|$ is less than or equal to
$$
|x|^{\frac{n}{2}}\sum_{\ga\in G(k)\cap C} |f_\infty(\ga \tilde{x}_{fin}^{-1}x_\infty)|.
$$
By the Artin-Whaples formula we have $|\tilde{x}_{fin}^{-1}|_\infty = |\tilde{x}_{fin}|_{fin}=|x_{fin}|_{fin}$ so that $|\tilde{x}_{fin}^{-1}x_\infty |_\infty =|x|$ and the sum is extended over a fixed lattice in $M(\A_\infty)$ which is intersected with $G(\A_\infty)$.
Since $f_\infty\in \CS(M(\A_\infty))$ the claim follows.
\qed

Let $\pi$ be an irreducible admissible representation of $G(\A)$ then $\pi =\otimes_v\pi_v$ and $\pi_v$ is of class one for almost all places $v$.
Let $S$ be a finite set of places such that  $\pi_v$ is of class one outside $S$.
Then $S$ contains all infinite places by definition.
Define the (partial) {\bf global $L$-function} of $\pi$ as
$$
L_S(\pi) \= \prod_{v\notin S}L(\pi_v),
$$
provided the product converges.

Now again $g\mapsto |g|^s$ for $s\in \C$ is admissible and so is
$\pi_s : g\mapsto |g|^s\pi(g)$.
Let $L(\pi,s) := L(\pi_s)$.

\begin{lemma}
Suppose $\pi$ is unitary, then for $\Re(s)>1$ the product $L_S(\pi,s)$ converges.
\end{lemma}

\prf
Since $\pi$ is unitary, all $\pi_v$ are, which implies that the Satake parameters $\chi_{\pi_v}(\varpi_j)$ have norm one.
\qed

Let $G(\A)^1$ be the kernel of the map $g\mapsto |g|$.
The exact sequence
$$
1\ra G(\A)^1)\begin{array}{c}i\\ \ra\\ {}\end{array} G(\A)\ra\R_+^\times\ra 1
$$
splits.
Moreover, there is splittings which make $G(\A)$ a direct product of $G(\A)^1$ and $\R_+^\times$.
For example fix a place $v_0 |\infty$ and let $s :\R_+^\times\ra G(\A)$ be given by $s(t)_v=\1$ if $v\ne v_0$ and $s(t)_{v_0} =t^{\rez{n}}\1$ then the map $(i,s) : G(\A)^1\times\R_+^\times \ra G(\A)$ is an isomorphism.
We will fix such a splitting from now on.
For $x\in G(\A)$ we will write $x=(x^1,|x|)$ in these coordinates.

The group $G(k)$ is a lattice in $G(\A)^1$.
A representation $\pi$ of $G(\A)^1$ is called {\bf automorphic} if it occurs as a subrepresentation in $L^2(G(k)\bs G(\A)^1)$.
With respect to our fixed splitting we are able to lift $\pi$ to a representation of $G(\A)$ by $\pi(x^1,|x|):=\pi(x^1)$ and we will thus consider it as  a representation of $G(\A)$.
Then $\pi$ decomposes into an infinite tensor product $\pi =\otimes_v\pi_v$.
We will write $Z$ for the image of the splitting map $s$ and we will thus identify $L^2(G(k)\bs G(\A)^1)$ to $L^2(ZG(k)\bs G(\A))$.
Write $R$ for the representation of $G(\A)$ on the latter.
That is, for $\ph\in L^2(ZG(k)\bs G(\A))$ and $y\in G(\A)$ we write $R(y)\ph(x) =\ph(xy)$.

Now assume $\pi$ is automorphic and fix a $G(\A)^1$-homomorphism $V_\pi\hookrightarrow L^2(G(k)\bs G(\A)^1)$.
Let $\ph\in L^2(ZG(k)\bs G(\A))$ be the image of a vector $\alpha =\otimes_v\alpha_v\in V_\pi =\otimes_vV_{\pi,v}$ such that $\alpha_v$ is a normalized class one vector at almost all places.
Further assume that $\ph$ is smooth and $\ph(1)\ne 0$.
The latter can be achieved by replacing $\ph(x)$ with $\ph(xy)$ for a suitable $y$ if necessary.

For any set of places $S$ let $\A_S$ be the restricted product of $k_v$, $v\in S$ and let $G_S=G(\A_S)$.
We consider $G_S$ as a subgroup of $G(\A)$.

\begin{lemma}\label{crucial}
Let $f\in\CS(M(\A))_0$ be of the form $f_S\otimes f^S$ for a finite set of places $S$, where $f_S\in \CS(M(\A_S))$ and $f^S =\prod_{v\notin S}\1_{M(\CO_v)}$.
Take $S$ so large that $\ph$ is of class one outside $S$, then for $\Re(s)>>0$ we have
$$
\int_{G(k)\bs G(\A)}E(f)(x)\ph_{\pi}(x)|x|^sd^*x \= L_S(\pi,s+\rez{2})\int_{G_S} f(x)\ph(x) |x|^{s+n/2}d^*x.
$$
The left hand side is entire in $s$. The second factor on the right hand side is entire if we assume that $f_S$ has compact support in $G_S$. 
The $L$-function is entire \cite{GJ}, hence the identity holds for all $s\in \C$.
\end{lemma}

\noindent
{\it Remark.}
The holomorphicity of the $L$-function is in \cite{GJ} only given for cuspidal ones.
The formula of the Lemma gives it for all automorphic $L$-functions.

\prf
For $\Re(s)>>0$ we compute that 
$$
\int_{G(k)\bs G(\A)}E(f)(x)\ph_{\pi}(x)|x|^sd^*x
$$
equals
\begin{eqnarray*}
\int_{G(\A)} f(x) |x|^{\frac{n}{2}+s}\ph(x)d^*x
	&=& \int_{G(\A)} f(x) |x|^{\frac{n}{2}+s}R(x)\ph(1)d^*x\\
	&=& \prod_{v\notin S}L(\pi_v,s+\rez{2})\int_{G_S}f(x)|x|^{\frac{n}{2}+s}\ph(x) d^*x.
\end{eqnarray*}
For the justification of this computation note that for $\Re(s)>>0$ the integral over $G(\A)$ converges absolutely.
\qed

In the notations of the lemma let
$$
\Delta_{S,\ph,s}(f) := \int_{G_S}f(x)\ph(x) |x|^{s+\frac{n}{2}}d^*x.
$$

For $\delta>0$ let $L^2_\delta(G(k)\bs G(\A))$ denote the space of meaurable functions $f$ on $G(k)\bs G(\A)$ with
$$
\int_{G(k)\bs G(\A)}|f(x)|^2\left(1+(\log |x|)^2\right)^{\delta /2}d^*x <\infty 
$$
modulo nullfunctions.
The sum $E$ defines a linear map $E:\CS(M(\A))_0\ra L^2_\delta(G(k)\bs G(\A))$
The group $G(\A)$ acts on both sides by right translation $R$ and it is easy to see that for any $y\in G(\A)$
$$
ER(y) \= |y|^{-n/2}R(y)E,
$$
so that the image of $E$ is an invariant subspace.

The pairing
$$
(f,g) \= \int_{G(k)\bs G(\A)}f(x)g(x) d^*x
$$
identifies $L_{-\delta}^2(G(k)\bs G(\A))$ to the dual space of $L_{\delta}(G(k)\bs G(\A))$.
The space $L_{\delta}^2(G(k)\bs G(\A))$ can be viewed as Hilbert space tensor product 
$$
L^2(G(k)\bs G(\A)^1)\otimes L^2_\delta(\R),
$$ 
where $L^2_\delta(\R)$ is the Fourier transform of the $\delta$-Sobolev space, i.e. the space of all functions $f$ on $\R$ with $\int_{\R}|f(x)|^2(1+x^2)^{\delta/2}dx <\infty$.

For an irreducible unitary representation $\pi$ of $G(\A)^1$ let 
$$
L^2(G(k)\bs G(\A)^1)(\pi)
$$ 
denote the isotypical component in $L^2(G(k)\bs G(\A)^1)$ and let 
$$
L^2_\delta(G(k)\bs G(\A))(\pi) \= L^2(G(k)\bs G(\A)^1)(\pi)\otimes L^2_\delta(\R)
$$
be its isotype in $L^2(G(k)\bs G(\A))$.
Let $\tilde{\pi}$ be the dual representation then the space $L^2_{-\delta}(G(k)\bs G(\A))(\tilde{\pi})$ is dual to $L^2_{\delta}(G(k)\bs G(\A))({\pi})$.
Let $\tilde{\ph} =\ph\otimes\psi\in L^2_{-\delta}(G(k)\bs G(\A))(\tilde{\pi})$,
then $\tilde{\ph}$ is orthogonal to the image of $E$ if and only if
$$
\int_{G(k)\bs G(\A)}E(f)(x)\ph(x)\psi(\log |x|)d^*x\= 0.
$$
We consider $\psi$ as a distribution and write formally $\psi(\log |x|)=\int_\R\hat{\psi}(t)|x|^{it}dt$.
The above becomes
$$
\int_{G(k)\bs G(\A)}\int_\R E(f)(x)\ph(x)\hat{\psi}(t)|x|^{it}dtd^*x
$$
which equals
$$
 \int_\R L_S(\pi,\rez{2}+it)\hat{\psi}(t)\Delta_{S,\ph,it}(f)dt
$$
if we plug in functions $f$ as in Lemma \ref{crucial}.
To justify this computation note that, since we are interested in the orthogonal space of an invariant space we can assume $R(h)\ph=\ph$ for some $h$ such that $\hat{h}$ has compact support.
Then $\hat{\psi}$ has compact support in $\R$.

At the finite places in $S$ let $f_v(x)=\1_{G(\CO_v)}\overline{\ph_v(x)}$.
Then there is a constant $c>0$ such that
$$
\Delta_{S,\ph ,it}(f) \= c\int_{G(\A_\infty)}f_\infty(x) \ph(x) |x|^{it+\frac{n}{2}}d^*x.
$$
We can choose $\ph$ so that there is a finite place $v$ with $R(G(\CO_v))\ph=0$ which implies that $f_\infty$ can be chosen arbitrarily and still $f\in\CS(M(\A))_0$.
If we let run $f_\infty$ through an approximate identity at the unit element we get that the distribution $L_S(\pi,\rez{2}+it)\hat{\psi}(t)$ on $\R$ is zero, i.e.
\begin{equation}\label{1}
L_S(\pi,\rez{2}+it)\hat{\psi}(t)\= 0.
\end{equation}
Now $\psi\in L^2_{-\delta}$ implies that its Fourier transform is a distribution of order $<\delta -1$.
For this recall that for $k\in\N$ the $k$-th derivative of the Dirac distribution $\delta^{(k)}(h)=h^{(k)}(0)$ has Fourier transform $(ix)^k$ which lies in $L_{-\delta}^2(\R)$ if and only if $k<\delta -1$.
Therefore equation (\ref{1}) is satisfied precisely for $\psi$ being a linear combination of $\delta^{(k)}(t)$ where $k<\delta -1$ and $\rez{2}+it$ is a zero of $L(\pi)$ of order $>k$.
We will interprete this as a spectral decomposition of the $t$-multiplication.
More precisely let $H_\pi$ be the Hausdorff cokernel of the composite map 
$$
\CS(M(\A))_0\begin{array}{c}E\\\ra\\{}\end{array} L_\delta^2(G(k)\bs G(\A))\begin{array}{c}Proj\\\ra\\{}\end{array}L_\delta^2(G(k)\bs G(\A))(\pi).
$$
Then $H_\pi$ is isomorphic to the orthogonal space of the image of $E$ in the dual space.
On $H_\pi$ we have the operator $D_\pi$ given by
$$
D_\pi\xi \= \lim_{\epsilon\ra 0} \rez{\epsilon}\left( R(s(e^\epsilon))-1\right)\xi,
$$
where $s :\R_+^\times\ra G(\A)$ is our fixed section.
The domain $\dom D_\pi$ is the set of $\xi$ for which the limit exists.
A computation shows that
$$
\norm{R(g)}_\delta \le 2^{\delta/4}(1+(\log |g^{-1}|)^2)^{\delta /4},
$$
so that $D_\pi$, which is the infinitesimal generator of $R(s(e^t))$, has purely imaginary spectrum.
Its resolvent $R_k = (D_\pi -k)^{-1}$ is given for $\Re(k)>0$ by
$$
R_k\= -\int_0^\infty R_\pi(s(e^t))e^{-k t}dt
$$
and for $\Re(k)<0$ by
$$
R_k =-\int_0^\infty R_\pi(s(e^{-t}))e^{k t}dt.
$$
The operator $D_\pi$ is closed since $v_n\ra v$ and $D_\pi v_n\ra y$ implies $v=R_k y-k R_k v$ lies in $\im R_k =\dom D_\pi$.
The operator $D_\pi$ commutes with the action of $G(\A)^1$.
Therefore its eigenspaces are $G(\A)^1$-modules, in fact, are multiples of $\pi$.
Let $m(\pi)$ be the multiplicity of $\pi$ in $L^2(G(k)\bs G(\A)^1)$.
Now recall that there is a definition \cite{GJ} for $L(\pi,s)$ in the global case by attaching factors at the ramified places in a way that the zero or poles of $L(\pi,s)$ along the line $\Re(s)=\rez{2}$ coincide with those of $L_S(\pi,s)$.
We have proven

\begin{theorem}
For $\delta >1$ the operator $D_\pi$ has discrete spectrum in $i\R$ consisting of all $i\rho\in i\R$ such that $\rez{2}+i\rho$ is a zero of $L(\pi,.)$. 
The eigenspace at $\rho$ is as a $G(\A)^1$-module a multiple of $\pi$ of multiplicity $m(\pi)n(\rho)$, where $n(\rho)$ is the largest integer $n<\delta -1$ such that $n<$ multiplicity of the zero $L(\pi ,\rez{2}+i\rho)$.
\end{theorem}
\qed

\tiny
\hspace{-20pt}
Math. Inst.\\ INF 288\\ 69120 Heidelberg\\ GERMANY

\end{document}